\documentclass[reqno,a4paper]{amsart}%
\usepackage{amssymb,url}
\usepackage[hypertex]{hyperref}
\allowdisplaybreaks[4]
\theoremstyle{plain}
\newtheorem{thm}{Theorem}
\theoremstyle{remark}
\newtheorem{rem}{Remark}
\DeclareMathOperator{\td}{d\mspace{-1mu}}
\date{}

\begin{document}

\title[Properties of extended remainder of Binet's formula]
{Some properties of extended remainder of Binet's first formula for logarithm of gamma function}

\author[F. Qi]{Feng Qi}
\address[F. Qi]{Department of Mathematics, College of Science, Tianjin Polytechnic University, Tianjin City, 300160, China}
\email{\href{mailto: F. Qi <qifeng618@gmail.com>}{qifeng618@gmail.com}, \href{mailto: F. Qi <qifeng618@hotmail.com>}{qifeng618@hotmail.com}, \href{mailto: F. Qi <qifeng618@qq.com>}{qifeng618@qq.com}}
\urladdr{\url{http://qifeng618.spaces.live.com}}

\author[B.-N. Guo]{Bai-Ni Guo}
\address[B.-N. Guo]{School of Mathematics and Informatics, Henan Polytechnic University, Jiaozuo City, Henan Province, 454010, China}
\email{\href{mailto: B.-N. Guo <bai.ni.guo@gmail.com>}{bai.ni.guo@gmail.com},
\href{mailto: B.-N. Guo <bai.ni.guo@hotmail.com>}{bai.ni.guo@hotmail.com}}
\urladdr{\url{http://guobaini.spaces.live.com}}

\begin{abstract}
In the paper, we extend Binet's first formula for the logarithm of the gamma function and investigate some properties, including inequalities, star-shaped and sub-additive properties and the complete monotonicity, of the extended remainder of Binet's first formula for the logarithm of the gamma function and related functions.
\end{abstract}



\thanks{The authors were supported in part by the China Scholarship Council and the Science Foundation of Tianjin Polytechnic University}

\thanks{This manuscript is a slightly revised version of the formally-published paper~\cite{Extended-Binet-remiander-comp.tex-Math-Slovaca}}


\maketitle

\section{Introduction}

\subsection{}
For positive numbers $x$ and $y$ with $y>x$, let
\begin{equation}
g_{x,y}(t)=\int_{x}^{y}u^{t-1}\td u=
\begin{cases}
\dfrac{y^{t}-x^{t}}{t}, & t\neq 0, \\
\ln y-\ln x, & t=0.
\end{cases}
\end{equation}
The reciprocal of $g_{x,y}(t)$ can be rewritten as
\begin{equation*}
\frac1{g_{e^a,e^b}(t)}=F_{a,b}(t)=
\begin{cases}
\dfrac{t}{e^{bt}-e^{at}},& t\ne 0,\\[0.6em]
\dfrac1{b-a},& t=0,
\end{cases}
\end{equation*}
where $a$ and $b$ are real numbers with $b>a$.
\par
It is well-known \cite[p.~11]{magnus} that Binet's first formula of $\ln\Gamma(x)$ for $x>0$ is given by
\begin{equation}\label{ebinet}
\ln \Gamma(x)= \left(x-\frac{1}{2}\right)\ln x-x+\ln \sqrt{2\pi}\,+\theta(x),
\end{equation}
where
\begin{equation*}
\Gamma(x)=\int^\infty_0t^{x-1} e^{-t}\td t
\end{equation*}
stands for Euler's gamma function and
\begin{equation}\label{ebinet1}
\theta(x)=\int_{0}^{\infty}\left(\frac{1}{e^{t}-1}-\frac{1}{t}
+\frac{1}{2}\right)\frac{e^{-xt}}{t}\td t
\end{equation}
for $x>0$ is called the remainder of Binet's first formula for the logarithm of the gamma function.
\par
In \cite{jmaa-ii-97, qx2}, some inequalities and completely monotonic properties of the function $g_{x,y}(t)$ were established and applied to construct Steffensen pairs in \cite{Gauchman-Steffensen-pairs}, to refine Gautschi-Kershaw's inequalities in \cite{notes-best.tex, notes-best.tex-rgmia}, and to study monotonic properties, logarithmic convexities and Schur-convexities of extended mean values $E(r,s;x,y)$ in \cite{schext, pams-62-rgmia, pams-62, schext-rgmia, ql, qx3}. See also \cite{cubo, cubo-rgmia} for related contents.
\par
In recent years, some inequalities and monotonic properties of the function
\begin{equation}\label{f-dfn-abs}
\frac{1}{t^2}-\frac{e^{-t}}{(1-e^{-t})^2}
\end{equation}
for $t>0$ and related ones were researched in \cite{best-constant-one-simple.tex, best-constant-one.tex, best-constant-one-simple-real.tex} and related references therein. These results were used in \cite{binet, remiander-Guo-Qi-Tamsui, psi-reminders.tex, psi-reminders.tex-rgmia} to consider completely monotonic properties of remainders of Binet's first formula, the psi function and related ones.
\par
Recently, logarithmic convexities of $g_{x,y}(t)$ and $F_{a,b}(t)$ were found in \cite{mon-element-exp-gen.tex, exp-funct-further.tex}. By virtue of these conclusions, some simple and elegant proofs for the logarithmic convexities, Schur-convexities of extended mean values $E(r,s;x,y)$ were simplified in \cite{emv-log-convex-simple.tex, Cheung-Qi-Rev.tex}.

\subsection{}
Now it is very natural to ask a question: Are there any relationship between the above studies? Direct computation yields
\begin{equation}\label{delta-a,b-(t)}
[g_{e^a,e^b}(t)]'=-[\ln F_{a,b}(t)]'=\frac{be^{bt}-ae^{at}}{e^{bt}-e^{at}}-\frac1t =\frac{b-a}{e^{(b-a)t}-1}-\frac1t+b\triangleq\delta_{a,b}(t)
\end{equation}
for $t\ne0$. Therefore, if taking $a=-\frac12$ and $b=\frac12$, then $\delta_{-1/2,1/2}(t)$ for $t>0$ equals the integrand in the remainder of Binet's first formula for the logarithm of the gamma function and the first order derivative of $\delta_{-1/2,1/2}(t)$ for $t>0$ also equals the function \eqref{f-dfn-abs}. These relationships connect closely the above three seemingly unrelated problems.
\par
If replacing $\delta_{-1/2,1/2}(t)$ by $\delta_{a,b}(t)$ on $t\in(0,\infty)$ for $b>a$
in \eqref{ebinet1}, then a more problem emerges: How to calculate the improper integral
\begin{equation}\label{int-binet-ext-1}
\int_0^\infty\biggl[\frac{b-a}{e^{(b-a)t}-1}-\frac1t+b\biggr]\frac{e^{-tx}}{t}\td t
\end{equation}
for $b>a$? The following Theorem~\ref{binet-first-ext-thm} answers this question affirmatively.

\begin{thm}\label{binet-first-ext-thm}
Let $b>a$, $\alpha>0$ and $x\in\mathbb{R}$ be real numbers. Then the improper integral \eqref{int-binet-ext-1} converges if and only if $a+b=0$ and
\begin{equation}\label{int-binet-ext-2}
\int_0^\infty\biggl(\frac{\alpha}{e^{\alpha t}-1}-\frac1t+\frac{\alpha}2\biggr)\frac{e^{-tx}}{t}\td t =\alpha\ln\Gamma\biggl(\frac{x}{\alpha}\biggr) -\biggl(x-\frac{\alpha}2\biggr)\ln\frac{x}{\alpha}+x-\frac{\alpha}2\ln(2\pi).
\end{equation}
\end{thm}

\begin{rem}
It is easy to see that the formula~\eqref{ebinet} is the special case $\alpha=1$ of \eqref{int-binet-ext-2}. So we call~\eqref{int-binet-ext-2} the extended remainder of Binet's first formula for the logarithm of the gamma function $\Gamma$.
\end{rem}

\subsection{}
It is well-known \cite{remiander-Guo-Qi-Tamsui, subadditive-qi-2.tex-arXiv, subadditive-qi-2.tex-rgmia} that a function $f$ is said to be completely monotonic on an interval $I$ if $f$ has derivatives of all orders on $I$ and
\begin{equation*}
(-1)^{n}f^{(n)}(x)\ge0
\end{equation*}
for $x \in I$ and $n \ge0$, that a function $f(x)$ is said to be star-shaped on $(0,\infty)$ if
\begin{equation*}
f(\alpha x)\le \alpha f(x)
\end{equation*}
for $x\in(0,\infty)$ and $0<\alpha<1$, that a function $f$ is said to be super-additive on $(0,\infty)$ if
\begin{equation*}
f(x+y)\ge f(x)+f(y)
\end{equation*}
for all $x,y>0$, and that a function $f$ is said to be sub-additive if $-f$ is super-additive.
\par
The function $\delta_{a,b}(t)$ defined by~\eqref{delta-a,b-(t)} has the following properties.

\begin{thm}\label{binet-first-ext-convex-thm}
Let $a$ and $b$ be real numbers with $a\ne b$, and let $0<\tau<1$.
\begin{enumerate}
\item
The function $\delta_{a,b}(t)$ is increasingly concave on $(0,\infty)$ and convex on $(-\infty,0)$;
\item
The inequality
\begin{equation}\label{delta-ab-ineq-1}
\delta_{a,b}(\tau t)<\delta_{a,b}(t)
\end{equation}
is valid on $(0,\infty)$; if either $a+b\ne0$ or $ab\ne0$, then the inequality
\eqref{delta-ab-ineq-1} is sharp;
\item
If $a+b\ge0$, then
\begin{equation}\label{delta-ab-ineq-2}
\tau\delta_{a,b}(t)<\delta_{a,b}(\tau t)
\end{equation}
is valid on $(0,\infty)$, i.e., the function $-\delta_{a,b}(t)$ is star-shaped; if $\max\{a,b\}\le0$, then the inequality \eqref{delta-ab-ineq-2} is reversed, i.e., the function $\delta_{a,b}(t)$ is star-shaped; if $a+b=0$, then the inequality \eqref{delta-ab-ineq-2} is sharp.
\end{enumerate}
\end{thm}

\begin{rem}
Some properties of special cases of the function $\delta_{a,b}(t)$ and related ones have been investigated and applied extensively in \cite{best-constant-one-simple.tex, emv-log-convex-simple.tex, mon-element-exp-gen.tex, remiander-Guo-Qi-Tamsui, best-constant-one.tex, exp-funct-further.tex, jmaa-ii-97, qx2, best-constant-one-simple-real.tex} and related references therein.
\end{rem}

\subsection{}
If denoting the extended remainder of Binet's first formula for the logarithm of the gamma function $\Gamma$ by
\begin{equation}\label{reminder-extended}
\theta_\alpha(x)=\int_0^\infty\biggl(\frac{\alpha}{e^{\alpha t}-1}-\frac1t+\frac{\alpha}2\biggr)\frac{e^{-tx}}{t}\td t
\end{equation}
for $\alpha>0$ and $x>0$, then formula \eqref{int-binet-ext-2} in Theorem~\ref{binet-first-ext-thm} can be simplified as
\begin{equation*}
\theta_\alpha(x)=\alpha\theta\biggl(\frac{x}\alpha\biggr)\quad \text{or}\quad
\theta_\alpha(\alpha x)=\alpha\theta(x).
\end{equation*}
This motivates us to study properties of $\theta_\alpha(x)$ and the function
\begin{equation*}
f_{p,q;\alpha}(x)=\theta_\alpha(px)-q\theta_\alpha(x)
\end{equation*}
on $(0,\infty)$, where $p>0$, $\alpha>0$ and $q\in\mathbb{R}$, which may be concluded as the following theorem.

\begin{thm}\label{extended-reminder-thm}
The extended remainder $\theta_\alpha(x)$ of Binet's first formula for the logarithm of the gamma function $\Gamma$ satisfies
\begin{equation}\label{theta-ineq}
\frac{(-1)^k}{(1+\lambda)^k}\theta_{\alpha}^{(k)}\biggl(\frac{x}{1+\lambda}\biggr) >\frac{(-1)^k}2\biggl[\frac1{\lambda^k}\theta_{\alpha}^{(k)}\biggl(\frac{x}\lambda\biggr) +\theta_\alpha^{(k)}(x)\biggr]
\end{equation}
for $x>0$, $\lambda>0$ with $\lambda\ne1$, $k\ge0$ and $\alpha>0$.
\par
The function $f_{p,q;\alpha}(x)$ is completely monotonic on $(0,\infty)$ if either $0<p\le1$ and $q\le 1$ or $p>1$ and $q\le\frac1p$; the function $-f_{p,q;\alpha}(x)$ is completely monotonic on $(0,\infty)$ if $p\ge1$ and $q\ge1$.
\par
The function $-\theta_{\alpha}(x)$ is star-shaped and $\theta_{\alpha}(x)$ is sub-additive.
\end{thm}

\begin{rem}
If taking $a=-\frac12$ and $b=\frac12$, then the results obtained in \cite{binet, Chen-Qi-Srivastava-09.tex, remiander-Guo-Qi-Tamsui, psi-reminders.tex, psi-reminders.tex-rgmia} can be deduced directly from Theorem~\ref{extended-reminder-thm}.
\end{rem}

\section{Proofs of theorems}

\begin{proof}[Proof of Theorem~\ref{binet-first-ext-thm}]
By transformations of integral variables, it easily follows that
\begin{align}\label{calculation}
&\quad\int_\varepsilon^\infty\biggl[\frac{b-a}{e^{(b-a)t}-1}-\frac1t+b\biggr]\frac{e^{-tx}}{t}\td t\notag\\
&=\int_{(b-a)\varepsilon}^\infty\biggl(\frac{b-a}{e^u-1}-\frac{b-a}u+b\biggr)\frac{e^{-ux/(b-a)}}u\td u\notag\\
&=(b-a)\int_{(b-a)\varepsilon}^\infty\biggl(\frac1{e^u-1}-\frac1u+\frac12\biggr)\frac{e^{-ux/(b-a)}}u\td u\notag\\
&\quad +\frac{a+b}2\int_\varepsilon^\infty\frac{e^{-ux/(b-a)}}u\td u\notag\\
\begin{split}
&=(b-a)\int_{(b-a)\varepsilon}^\infty\biggl(\frac1{e^u-1}-\frac1u+\frac12\biggr)\frac{e^{-ux/(b-a)}}u\td u\\
&\quad+\frac{a+b}2\int_{(b-a)\varepsilon}^\infty\frac{e^{-ux}}u\td u,
\end{split}
\end{align}
where $\varepsilon>0$.
\par
By virtue of Binet's first formula for $\ln\Gamma(z)$ in \eqref{ebinet1}, the integral in the first term of \eqref{calculation} may be calculated as
\begin{equation}\label{binet's-first-rev}
\begin{aligned}
&\quad\lim_{\varepsilon\to0^+}\int_{(b-a)\varepsilon}^\infty\biggl(\frac1{e^u-1}-\frac1u+\frac12\biggr) \frac{e^{-ux/(b-a)}}u\td u \\
&=\int_0^\infty\biggl(\frac1{e^u-1}-\frac1u+\frac12\biggr) \frac{e^{-ux/(b-a)}}u\td u \\
&=\ln\Gamma\biggl(\frac{x}{b-a}\biggr)-\biggl(\frac{x}{b-a}-\frac12\biggr) \ln\frac{x}{b-a}+\frac{x}{b-a}-\frac12\ln(2\pi).
\end{aligned}
\end{equation}
Furthermore, the second integral in~\eqref{calculation} satisfies
\begin{align*}
\lim_{\varepsilon\to0^+}\int_{(b-a)\varepsilon}^\infty\frac{e^{-ux}}u\td u
=\lim_{\varepsilon\to0^+}\int_{(b-a)x\varepsilon}^\infty t^{-1}e^{-t}\td t
=\int_0^\infty t^{-1}e^{-t}\td t
\end{align*}
which is divergent. As a result, the improper integral \eqref{int-binet-ext-1} is
convergent if and only if $a+b=0$.
\par
Taking $a=-b$ in \eqref{int-binet-ext-1} and \eqref{binet's-first-rev} and simplifying yields formula \eqref{int-binet-ext-2}. The proof of Theorem~\ref{binet-first-ext-thm} is complete.
\end{proof}

\begin{proof}[Proof of Theorem~\ref{binet-first-ext-convex-thm}]
Straightforward computation gives
\begin{gather*}
\begin{split}
\delta_{a,b}'(t)&=\frac1{t^2}-\frac{(a-b)^2e^{(a+b)t}}{(e^{at}-e^{bt})^2},\\
\delta_{a,b}''(t)&=\frac{(a-b)^3e^{(a+b)t} (e^{at}+e^{bt})}{(e^{a t}-e^{b t})^3}-\frac2{t^3}
\end{split}\\
=\frac{2e^{3(a+b)t/2}}{t^3}\biggl(\frac{at-bt}{e^{at}-e^{bt}}\biggr)^3 \biggl\{\frac{e^{(a-b)t/2}+e^{(b-a)t/2}}2 -\biggl[\frac{e^{(a-b)t/2}-e^{(b-a)t/2}}{(a-b)t}\biggr]^3\biggr\}\\
\triangleq\frac{2e^{3(a+b)t/2}}{t^3} \biggl(\frac{at-bt}{e^{at}-e^{bt}}\biggr)^3Q\biggl(\frac{a-b}2t\biggr).
\end{gather*}
Lazarevi\'c's inequality in \cite[p.~131]{bullen-dic} and \cite[p.~300]{3rded} tells us that
\begin{equation*}
Q(t)=\frac{e^{-t}+e^t}{2} -\biggl(\frac{e^t-e^{-t}}{2t}\biggr)^3
=\cosh t -\biggl(\frac{\sinh t}{t}\biggr)^3<0
\end{equation*}
for $t\in\mathbb{R}$ with $t\ne0$. Hence $\delta_{a,b}''(t)<0$ on $(0,\infty)$ and
$\delta_{a,b}''(t)>0$ on $(-\infty,0)$. The convexity and concavity of
$\delta_{a,b}(t)$ are proved.
\par
Since $\delta_{a,b}''(t)<0$ on $(0,\infty)$, the derivative $\delta_{a,b}'(t)$ is
decreasing on $(0,\infty)$ for all real numbers $a$ and $b$ with $a\ne b$. Since
\begin{equation*}
\delta_{a,b}'(t)=\frac1{t^2}-\frac{(a-b)^2e^{(b-a)t}}{[1-e^{(b-a)t}]^2} =\frac1{t^2}-\frac{(a-b)^2e^{(a-b)t}}{[1-e^{(a-b)t}]^2},
\end{equation*}
it is easy to obtain that
$$
\lim_{t\to\infty}\delta_{a,b}'(t)=0.
$$
Consequently, the function $\delta_{a,b}'(t)$ is positive, and so $\delta_{a,b}(t)$ is
increasing, on $(0,\infty)$. This means that inequality \eqref{delta-ab-ineq-1} holds
for $0<\tau<1$ and $t>0$.
\par
From
\begin{equation*}
\delta_{a,b}(t)=\frac{be^{bt}-ae^{at}}{e^{bt}-e^{at}}-\frac1t =\frac{be^{(b-a)t}-a}{e^{(b-a)t}-1}-\frac1t =\frac{b-ae^{(a-b)t}}{1-e^{(a-b)t}}-\frac1t,
\end{equation*}
it follows easily that
\begin{equation*}
\lim_{t\to\infty}\delta_{a,b}(t)=\max\{a,b\}.
\end{equation*}
L'H\^opital's rule gives
\begin{equation*}
\begin{split}
\lim_{t\to0^+}\delta_{a,b}(t) &=\lim_{t\to0^+}\frac{t(be^{bt}-ae^{at})-e^{bt}+e^{at}}{t(e^{bt}-e^{at})} \\ &=\lim_{t\to0^+}\frac{b^2e^{bt}-a^2e^{at}} {(be^{bt}-ae^{at})+(e^{bt}-e^{at})/t}\\
&=\frac{a+b}2.
\end{split}
\end{equation*}
For $0<\tau<1$, let
$$
h_{a,b}(t)=\delta_{a,b}(\tau t)-\tau \delta_{a,b}(t)
$$
for $t>0$. It is obvious that
\begin{equation*}
\lim_{t\to0^+}h_{a,b}(t)=\frac{(1-\tau)(a+b)}2\quad \text{and}\quad \lim_{t\to\infty}h_{a,b}(t)=(1-\tau)\max\{a,b\}.
\end{equation*}
Since $\delta_{a,b}'(t)$ is decreasing, then
$$
h'_{a,b}(t)=\tau\bigl[\delta_{a,b}'(\tau t)-\delta_{a,b}'(t)\bigr]>0,
$$
and so $h_{a,b}(t)$ is strictly increasing. If $a+b\ge0$, then $h_{a,b}(t)>0$ on $(0,\infty)$ and inequality \eqref{delta-ab-ineq-2} is valid. If $\max\{a,b\}\le0$, then inequality \eqref{delta-ab-ineq-2} is reversed.
\par
It is apparent that
$$
\lim_{t\to \infty}\frac{\delta_{a,b}(\tau t)}{\delta_{a,b}(t)}=1
$$
if $ab\ne0$, which implies that inequality \eqref{delta-ab-ineq-1} is sharp. If
$a+b\ne0$, then it is clear that
$$
\lim_{t\to0^+}\frac{\delta_{a,b}(\tau t)}{\delta_{a,b}(t)}=1,
$$
which also implies that inequality \eqref{delta-ab-ineq-1} is sharp.
\par
By L'H\^opital's rule, it is not difficult to obtain that
$$
\lim_{t\to0^+}\delta_{a,b}'(t)=\frac{(a-b)^2}{12}.
$$
If $a+b=0$ and $a\ne b$, then
$$
\lim_{t\to0^+}\frac{\delta_{a,b}(\tau t)}{\delta_{a,b}(t)}
=\lim_{t\to0^+} \frac{\tau\delta_{a,b}'(\tau t)}{\delta_{a,b}'(t)}=\tau,
$$
which means that inequality \eqref{delta-ab-ineq-2} is sharp. The proof of
Theorem~\ref{binet-first-ext-convex-thm} is complete.
\end{proof}

\begin{proof}[Proof of Theorem~\ref{extended-reminder-thm}]
From the concavity of $\delta_{a,b}(t)$ on $(0,\infty)$, it follows that
\begin{equation*}
\delta_{-\alpha/2,\alpha/2}\biggl(\frac{(1+\lambda)t}2\biggr)
>\frac{\delta_{-\alpha/2,\alpha/2}(\lambda t)+\delta_{-\alpha/2,\alpha/2}(t)}2
\end{equation*}
for $t>0$ and positive numbers $\alpha$ and $\lambda\ne1$. Multiplying by the factor $t^{k-1}e^{-tx}$ for any nonnegative integer $k\ge0$ and integrating from $0$ to $\infty$ on both sides of the above inequality yields
\begin{multline*}
\int_0^\infty\biggl[\frac{\alpha}{e^{\alpha(1+\lambda)t/2}-1} -\frac2{(1+\lambda)t}+\frac{\alpha}2\biggr]t^{k-1}e^{-tx}\td t\\ >\frac12\biggl[\int_0^\infty\biggl(\frac{\alpha}{e^{\alpha\lambda t}-1} -\frac1{\lambda t} +\frac{\alpha}2\biggr)t^{k-1}e^{-tx}\td t\\
+\int_0^\infty\biggl(\frac{\alpha}{e^{\alpha t}-1}
-\frac1t+\frac{\alpha}2\biggr)t^{k-1}e^{-tx}\td t\biggr]
\end{multline*}
which can be rewritten by transformations of integral variables as
\begin{multline*}
\int_0^\infty\biggl(\frac{\alpha}{e^{\alpha u}-1} -\frac1u+\frac{\alpha}2\biggr)\frac{u^{k}} {(1+\lambda)^k}\cdot\frac{e^{-xu/(1+\lambda)}}u\td u\\ >\frac12\biggl[\int_0^\infty\biggl(\frac{\alpha}{e^{\alpha u}-1} -\frac1u+\frac{\alpha}2\biggr)\frac{u^{k}} {\lambda^k}\cdot\frac{e^{-xu/\lambda}}u\td u\\
+\int_0^\infty\biggl(\frac{\alpha}{e^{\alpha t}-1}
-\frac1t+\frac{\alpha}2\biggr)t^{k}\frac{e^{-tx}}t\td t\biggr].
\end{multline*}
Substituting formula \eqref{reminder-extended} and its derivatives into the above inequalities leads to
\begin{equation*}
(-1)^k\biggl[\theta_{\alpha}\biggl(\frac{x}{1+\lambda}\biggr)\biggr]^{(k)} >\frac12\biggl\{(-1)^k\biggl[\theta_{\alpha}\biggl(\frac{x}\lambda\biggr)\biggr]^{(k)} +(-1)^k\theta_\alpha^{(k)}(x)\biggr\}.
\end{equation*}
As a result, inequalities in \eqref{theta-ineq} follow.
\par
Easy calculation yields
\begin{align*}
f_{p,q;\alpha}(x)&=\int_{0}^{\infty}\delta_{-\alpha/2,\alpha/2}(t) \frac{e^{-pxt}}{t}\td t -q\int_{0}^{\infty}\delta_{-\alpha/2,\alpha/2}(t) \frac{e^{-xt}}{t}\td t\\
&=\int_{0}^{\infty}\delta_{-\alpha/2,\alpha/2}\biggl(\frac{t}p\biggr) \frac{e^{-xt}}{t}\td t-q\int_{0}^{\infty}\delta_{-\alpha/2,\alpha/2}(t) \frac{e^{-xt}}{t}\td t\\
&=\int_{0}^{\infty}\biggl[\delta_{-\alpha/2,\alpha/2}\biggl(\frac{t}p\biggr) -q\delta_{-\alpha/2,\alpha/2}(t)\biggr] \frac{e^{-xt}}{t}\td t\\
&\triangleq\int_{0}^{\infty}h_{p,q;\alpha}(t)\frac{e^{-xt}}{t}\td t.
\end{align*}
By virtue of properties of $\delta_{a,b}(t)$ obtained in
Theorem~\ref{binet-first-ext-convex-thm}, it follows by standard arguments that
\begin{enumerate}
\item
$h_{p,q;\alpha}(t)\ge0$ if either $0<p\le1$ and $q\le 1$, or $p>1$ and $q\le\frac1p$, or $0<q<1$ and $q\le\frac1p$;
\item
$h_{p,q;\alpha}(t)\le0$ if $p\ge1$ and $q\ge1$.
\end{enumerate}
It is clear that
\begin{enumerate}
\item
if $h_{p,q;\alpha}(t)\ge0$ then $f_{p,q;\alpha}(x)$ is completely monotonic on
$(0,\infty)$;
\item
if $h_{p,q;\alpha}(t)\le0$ then $-f_{p,q;\alpha}(x)$ is completely monotonic on
$(0,\infty)$.
\end{enumerate}
As a result, the completely monotonic properties of $f_{p,q;\alpha}(x)$ is proved.
\par
It is easy to see that the star-shaped properties of the function $\theta_\alpha(x)$
follow from those of the function $\delta_{-\alpha/2,\alpha/2}(t)$ and formula
\eqref{reminder-extended}.
\par
In \cite[p. 453]{mo}, it was presented that a star-shaped function must be
super-additive, therefore the function $\theta_\alpha(x)$ is also sub-additive. The
proof of Theorem~\ref{extended-reminder-thm} is complete.
\end{proof}

\begin{rem}
Dividing both sides of \eqref{int-binet-ext-2} by $\alpha>0$ gives
\begin{equation}\label{int-binet-ext-2-divid}
\int_0^\infty\biggl(\frac1{e^{\alpha t}-1}-\frac1{\alpha t}+\frac12\biggr)\frac{e^{-tx}}{t}\td t =\ln\Gamma\biggl(\frac{x}{\alpha}\biggr) -\biggl(\frac{x}\alpha-\frac12\biggr)\ln\frac{x}{\alpha}+\frac{x}\alpha-\frac12\ln(2\pi).
\end{equation}
Further taking $x=\alpha y$ in \eqref{int-binet-ext-2-divid} yields
\begin{equation}\label{int-binet-ext-2-y}
\int_0^\infty\biggl(\frac1{e^{\alpha t}-1}-\frac1{\alpha t}+\frac12\biggr)\frac{e^{-\alpha ty}}{t}\td t =\ln\Gamma(y) -\biggl(y-\frac12\biggr)\ln y+y-\frac12\ln(2\pi).
\end{equation}
Since
\begin{equation}
\int_0^\infty\biggl(\frac1{e^{\alpha t}-1}-\frac1{\alpha t}+\frac12\biggr)\frac{e^{-\alpha ty}}{t}\td t =\int_0^\infty\biggl(\frac1{e^{t}-1}-\frac1{t}+\frac12\biggr)\frac{e^{-ty}}{t}\td t,
\end{equation}
so the identity~\eqref{int-binet-ext-2} is essentially equivalent to~\eqref{ebinet}.
\end{rem}
%

\end{document}